*Original Article*

# Matroid, Ideal, Ultrafilter, Tangle, and so on :Reconsideration of Obstruction to Linear Decomposition

Takaaki Fujita

*Graduate school of Science and Technology, Gunma University, 1-5-1 Tenjin-cho Kiryu Gunma, Japan.*

*Corresponding Author : t171d603@gunma-u.ac.jp*

***Abstract*** *- The investigation of width parameters in both graph and algebraic contexts has attracted considerable interest. Among these parameters, the linear branch width has emerged as a crucial measure. In this concise paper, we explore the concept of linear decomposition, specifically focusing on the single filter in a connectivity system. Additionally, we examine the relevance of matroids, antimatroids, and greedoids in the context of connectivity systems. Our primary objective in this study is to shed light on the impediments to linear decomposition from multiple perspectives.*

## 1. Introduction

Graph theory explores mathematical structures called graphs, which model relationships between pairs of objects using vertices (nodes) and edges (links). This field analyzes graphs using parameters like degree, diameter, and chromatic number, which quantify various graph aspects for both theoretical study and practical applications.

Recently, the study of width parameters has garnered significant attention in graph theory and algebraic contexts, as evidenced by extensive research [1- 9, 12-14,17, 21, 22, 28-34, 48-50,56,57,60,61]. Width parameters are derived from tree-like structures, known as decompositions, on graphs. For instance, tree width, a well-known width parameter, measures how closely a graph resembles a tree and finds applications across databases, networks, artificial intelligence, protein structures, and programming languages [36-41]. Another crucial width parameter, linear branch width, helps characterize the complexity of mathematical objects such as graphs and matroids [31-33, 51].

The concept of a "Linear tangle" was introduced by Bienstock et al. [10] as part of a strategy for the game of cops and robbers. This concept and its related strategies reveal that an order $k + 1$ linear tangle in a connectivity system $(X,f)$, or a system with linear branch width $k$, can be derived using graph exploration games. A connectivity system is defined by a finite set X and a symmetric submodular function f [16]. "Tangle," a concept introduced by Robertson et al., exhibits properties similar to Linear tangle and is dual to the "branch width" width parameter [11]. Thus, the exploration of Tangle and Linear tangle, both closely linked to width parameters, remains a crucial area of research.

"Filter" is a well-known concept in the realms of topology and algebra. In simple terms, a filter can be interpreted as a collection of sets that contain a particular element, and it serves as a convenient tool for discussing convergence properties in mathematics. In the domain of Boolean algebra, filters that are maximal are referred to as "ultrafilters". It is established that ultrafilters on a connectivity system $(X,f)$ with conditions of a symmetric submodular function have a profound connection with Tangles and branch width [15]. It is both surprising and highly intriguing that by imposing constraints of symmetric submodular functions on algebraic concepts such as filters, one can obtain dual theorems for width parameters.

The purpose of this concise manuscript is to emphasize the practical significance of studying linear decomposition. To accomplish this, we introduce and define linear concepts that are closely associated with linear tangle and linear width, and we establish their equivalence. Specifically, we delve into the connection between filters and linear width. Within this exploration, we integrate an algorithmic concept known as "single element deletion" into filters. Single element deletion involves removing one element from a set while preserving its original properties. It is a fundamental and indispensable operation in the realm of matroid theory, and can be interpreted as a crucial aspect of algorithm development within this field. Throughout this manuscript, we refer to filters that incorporate single element deletion as "single ultrafilter." we demonstrate the dual relationship between single ultrafilters and linear branch-width.

Additionally, we consider about matroid, antimatroid, greedoid on connectivity system. Matroids are a discrete representation of linear algebra concepts such as linearly independent and linearly dependent sets. These concepts also have deep relation to linear width.

Let me explain the structure of this paper. In Section 2, we provide necessary definitions for this study. In Section 3, we characterize single filters. Section 4 examines the relationship between linear tangles and single ultrafilters. Section 5 discusses Matroids, Antimatroids, and Greedoids. Finally, in Section 6, we present the main theorem of this paper.

## 2. Definitions and Notations in this Paper

This section provides mathematical definitions of each concept. Here, I will describe the basic notation used in mathematics. In this paper, a set is defined as a collection of distinct elements or objects, treated as a single entity and typically denoted by curly braces. A subset is a set whose elements are all contained within another set. Boolean algebra $(X, \cup, \cap)$ represents a mathematical structure comprised of a set $X$, along with operations of union $(\cup)$ and intersection $(\cap)$, which satisfy specific axioms. Furthermore, this paper focuses on finite sets.

In this paper, we use expressions like $A \subseteq X$ to indicate that $A$ is a subset of $X$, $A \cup B$ to represent the union of two subsets $A$ and $B$ (both of which are subsets of $X$), and $A = \emptyset$ to signify an empty set. Specifically, $A \cap B$ denotes the intersection of subsets $A$ and $B$. Similarly, $A \setminus B$ represents the difference between subsets $A$ and $B$. The powerset of a set $A$, denoted as $2^A$, is the set of all possible subsets of $A$, including the empty set and $A$ itself. A singleton set, also known as a single-element set or a unit set, is a set containing exactly one element. No more, no less. In mathematical notation, if the element is $e$, then the singleton set is denoted as $\{e\}$. The notation $|A|$ represents the cardinality of set $A$, indicating the number of elements contained within $A$.

### *2.1. Filters on Boolean Algebra*

The definition of a filter in a Boolean algebra $(X, \cup, \cap)$ is given below.

**Definition 1:** In a Boolean algebra $(X, \cup, \cap)$, a set family $F \subseteq 2^X$ satisfying the following conditions is called a filter on the carrier set X.
(FB1) $A, B \in F \Rightarrow A \cap B \in F$,
(FB2) $A \in F, A \subseteq B \subseteq X \Rightarrow B \in F$,
(FB3) $\emptyset$ is not belong to $F$.

In a Boolean algebra $(X, \cup, \cap)$, A maximal filter is called an ultrafilter and satisfies the following axiom (FB4):
(FB4) $\forall A \subseteq X$, *either* $A \in F$ *or* $X / A \in F$.

And filter (ultrafilter) is non-principal if the filter (ultrafilter) satisfies following axiom:
(FB5) $A \notin F$ for all $A \subseteq X$ with $|A| = 1$.

Ultrafilters have wide-ranging and significant applications in a variety of fields, including topology, algebra, logic, set theory, lattice theory, matroid theory, graph theory, combinatorics, measure theory, model theory, and functional analysis (see [23-27]). They are extensively employed to understand mathematical structures, explore infinite objects, and provide powerful tools for reasoning and analysis in various domains.

An equivalent form of a given ultrafilter $U \subseteq 2^X$ is clearly manifested as a two-valued morphism. Note that a two-valued morphism is a function that maps elements to only two possible values, typically $0$ and $1$, used in decision processes or to characterize binary properties in set systems. This relationship is defined through a function $m$ where $m(A) = 1$ if A is an element of $U$, and $m(A) = 0$ otherwise.

### *2.2. Symmetric Submodular Function*

A symmetric submodular function is a set function that decreases with set addition and is invariant under set complement, useful in optimization and graph theory.Below is the definition of a symmetric submodular function. It is worth noting that while symmetric submodular functions can typically have real values, this paper specifically concentrates on the subset of functions that exclusively deal with natural numbers.

**Definition 2:** Let $X$ be a finite set. A function $f: X \rightarrow \mathbb{N}$ is called symmetric submodular if it satisfies the following conditions:
· $\forall A \subseteq X, f(A) = f(X/A)$.
· $\forall A, B \subseteq X, f(A) + f(B) \geq f(A \cap B) + f(A \cup B)$.

It is known that a symmetric submodular function $f$ satisfies the following properties. This lemma will be utilized in the proofs of lemmas and theorems presented in this paper.

**Lemma 3(cf.[16]) :** A symmetric submodular function $f$ satisfies:
1. $\forall A \subseteq X, f(A) \geq f(\emptyset) = f(X)$.
2. $\forall A, B \subseteq X, f(A) + f(B) \geq f(A \backslash B) + f(B \backslash A)$.

**Proof.** From the properties of a symmetric submodular function, the following results can be obtained.
(1): $f(A) +f(A) = f(A) +f(A) \geq f(A \cup A) +f(A \cap A) = f(X) +f(\emptyset) =f(\emptyset)+f(\emptyset)$.
(2): $f(A)+f(B) = f(A)+f(B) \geq f(A \cup B)+f(A\cap B) = f(B\backslash A)+f(A\backslash B) = f(B\backslash A) + f(A\backslash B)$. This proof is completed.

The exploration of symmetric submodular functions has made a significant contribution to the advancement of graph algorithms and our understanding of the structure and properties of submodular functions as a whole.
To facilitate comprehension, let's present a familiar example of a symmetric submodular function.

**Example 4:** Let's examine a graph $G = (V, E)$ and define the function $f_{\Phi}(A)$ for each $A \subset E$. This function calculates the number of vertices $v \in V$ that are adjacent to both edges belonging to set $A$ and edges not belonging to $A$. In simpler terms, $f_{\Phi}(A)$ counts the vertices that share a connection with both types of edges. It is evident that the function $f_{\Phi}(A)$ qualifies as a symmetric submodular function.

In this short paper, a pair $(X, f)$ of a finite set $X$ and a symmetric submodular function $f$ is called a connectivity system (cf. [16]). In this paper, we use the notation $f$ for a symmetric submodular function, a finite set $X$, and a natural number $k$. A set $X$ is $k$-efficient if $f(X) \leq k$.

### *2.3. Linear Tangle and Linear Decomposition*

The definition of a linear tangle is given below.

**Definition 5 [10]:** Let $X$ represent a finite set and $f$ denote a symmetric submodular function. A linear tangle of order $k+1$ on a connectivity system $(X,f)$ is a family $L$ of $k$-efficient subsets of $X$, satisfying the following axioms:
(L1) $\emptyset \in L$
(L2) For each $k$-efficient subset $A$ of $X$, exactly one of $A$ or $X/A$ in $L$
(L3) If $A,B \in L$, $e \in X$, and $f(\{e\}) \leq k$, then $A \cup B \cup \{e\} \neq X$ holds.

The concept mentioned in literature [10] pertains to discussions on graphs. Specifically, it asserts that a linear tangle of order $k + 1$ is also considered a linear tangle of order $k$. Moreover, it should be noted that if $k$ exceeds a certain threshold, a linear tangle cannot exist.

Linear tangle is a linearly restricted version of tangle. In literature [35], the concept of Tangle kit is defined, and it is known that the following holds true. Naturally, this also holds true when considering the linear world. In this paper, we will refer to the linearly restricted Tangle kit as the Linear Tangle kit.

**Theorem 6[35].** Let $X$ be a finite set and $f$ be a symmetric submodular function. There exists a tangle of order $k+1$ if and only if there exists a tangle kit of order $k+1$.

Let's present the definitions of linear decomposition

**Definition 7:** Let $C$ be a caterpillar, which is a tree with interior vertices of degree $3$ and leaves of degree $1$. Let $C$ be the path $(l_1, b_2, b_3, . . . , b_{n-1}, l_n)$. For $2 \leq i \leq n-1$, the subgraph of $C$ induced by $\{b_{i-1}, b_i, b_{i+1}\}$ is a connectivity system $(X, f)$. The Linear Decomposition of $C$ is a caterpillar that partitions the elements of $X$ into sets $\{e_1\}, \{e_2\}, ..., \{e_{n-1}\}, \{e_n\}$ as follows: for each $1 \leq i \leq n-1$, let $wi := f(\{e_1, ..., e_i\})$. The width of the Linear Decomposition is defined as $max\{w_1, ..., w_{n-1}, f(e_1), ..., f(e_{n-1}), f(e_n)\}$, and the linear width of $(X, f)$ is the smallest width among all Linear Decompositions of $(X, f)$.

Linear decomposition is a technique for partitioning a graph into simpler structured. As we said in the introduction section, this technique has been used to solve several long-standing problems in graph theory, including the Graph Minors Project by Robertson and Seymour. The study of linear decomposition continues to be an active area of research in graph theory, with applications in algorithms, complexity theory, and structural graph theory.

The following theorems and lemmas regarding linear decomposition/linear tangle are already known:

**Theorem 8 [10]:** If $X$ is a finite set and $f$ is a symmetric submodular function, linear-width of the connectivity system $(X, f)$ is at most $k$ if and only if no linear tangle of order $k+1$ exists.

**Theorem 9[31]:** Let $X$ be a finite set and $f$ be a symmetric submodular function. Under the assumption that $f(\{e\}) \leq k$ for every $e \in X$, $F$ is a linear tangle of order $k+1$ on $(X,f)$ iff F is a linear obstacle of order $k+1$ on $(X,f)$.

**Theorem 10[13]:** Let $X$ be a finite set and $f$ be a symmetric submodular function. Under the assumption that $f(\{e\}) \leq k$ for every $e \in X$, $F$ is a linear tangle of order $k+1$ on $(X,f)$ iff F is a maximal single ideal of order $k+1$ on $(X,f)$.

**Theorem 11[32]:** If $X$ is a finite set and $f$ is a symmetric submodular function, linear tangle of order $k+1$ exists if and only if maximal linear loose tangle of order $k+1$ exists.

## 3. Characterization of Single Ultrafilter

We introduce a new concept of single ultrafilter on a connectivity system $(X, f)$ by following the definition of ultrafilters on Boolean algebras.

**Definition 13:** Let $X$ represent a finite set and $f$ denote a symmetric submodular function. A subset $S \subseteq 2^X$ is called an order $k + 1$ single filter on a connectivity system $(X, f)$ if it satisfies the following conditions:

(S1) For any $A \in S$, $e \in X$, if $f(\{e\}) \leq k$ and $f(A \cap (X/\{e\})) \leq k$, then $A \cap (X/\{e\}) \in S$.

(S2) For any $A \in S$ and $A \subset B \subseteq X$, if $f(B) \leq k$, then $B \in S$.

(F3) $\emptyset$ is not belong to $S$.

**Definition 14:** Let $X$ represent a finite set and $f$ denote a symmetric submodular function. A subset $S \subseteq 2^X$ is called an order $k + 1$ single ultrafilter on a connectivity system $(X, f)$ if it satisfies the following conditions:

(S1) For any $A \in S$, $e \in X$, if $f(\{e\}) \leq k$ and $f(A \cap (X/\{e\})) \leq k$, then $A \cap (X/\{e\}) \in S$.

(S2) For any $A \in S$ and $A \subset B \subseteq X$, if $f(B) \leq k$, then $B \in S$.

(F3) $\emptyset$ is not belong to $S$.

(S4) For any $A \subseteq X$, if $f(A) \leq k$, either $A \in S$ or $(X/A) \in S$.

And single filter (single ultrafilter) is non-principal if the single filter (single ultrafilter) satisfies following axiom:

(F5) $A \notin S$ for all $A \subseteq X$ with $|A| = 1$, $f(A) \leqq k$.

The definition of a single ultrafilter on a connectivity system $(X, f)$ differs from that of an ultrafilter on a Boolean algebra $(X, \cup, \cap)$ in two aspects:

- Each axiom includes the condition of being a symmetric submodular function.
- The set B in axiom (F1) is replaced with a singleton set. This replacement is known as "single-element-deletion." Single-element-deletion is a widely used concept in the world of matroids, along with the idea of single-element-extension (cf. [18-20]). It is fascinating to note the connection between seemingly unrelated concepts in Boolean algebra, graph theory, and the concept of matroids.

For reference, let's introduce the concept of single-element deletion in a filter on $(X, f)$ with an additional constraint of submodularity. Consider a subset $S \subseteq 2^X$ that forms an order $k + 1$ single ultrafilter on $(X, f)$. Single-element deletion, referred to as the linear filter's single-element deletion, is denoted as $S \setminus \{e\}$. It is defined as $\{A\{e\} \mid f(A\{e\}) \leq k, f(\{e\}) \leq k, A \in S\}$. This definition is a natural extension of the fundamental idea of deletion/single-element deletion in matroid theory to a connectivity system.

In this paper, we assume that $f(\{e\}) \leq k$ for any element e, allowing us to rewrite $S \setminus \{e\}$ as $\{A\{e\} \mid f(A\{e\}) \leq k, A \in S\}$. According to the axiom (S1), a linear filter $S$ possesses the capability to "filter out" elements $e$ from any set $A \in S$, provided that $f(\{e\}) \leq k$. Moreover, since every element $e \in X$ satisfies $f(\{e\}) \leq k$, we can continue applying the filtering process to remove element e as long as the resulting set remains $k$-efficient. This operation aligns perfectly with the concept of single-element deletion in matroid theory. In fact, by replacing the axiom (S1) with the following (SD1), we precisely capture the idea of single-element deletion on $(X, f)$:

(SD1) For $A \in S, e \in X,$ if $f(A\{e\}) \leq k,$ then $A \setminus \{e\} \in S.$

We examine the maximality of single filters and single ultrafilters. The following lemma demonstrates the maximality of single ultrafilters.

**Lemma 15:** Let $X$ represent a finite set and $f$ denote a symmetric submodular function. Under the assumption that $f(\{e\}) \leq k$ for every $e \in X$, If subset $S \subseteq 2^X$ is order $k + 1$ single ultrafilter on $(X, f)$, then subset $S \subseteq 2^X$ is maximal single filter of order $k + 1$ on $(X, f)$.

**Proof:** Let $S$ be an order $k + 1$ single ultrafilter on $(X, f)$. We need to show that S is a maximal single filter of order $k + 1$ on $(X, f)$.

Suppose, for the sake of contradiction, that there exists a single filter $F$ such that $S \subset F$ and $F$ is of order $k + 1$. We will show that this leads to a contradiction, thus proving that $S$ is maximal.

First, note that F must satisfy conditions (S1) and (S2) of Definition, since it is a single filter of order $k + 1$. We will show that F also satisfies condition (S4) of Definition.

Suppose there exists a set $A \subseteq X$ such that $f(A) \leq k$ and neither $A \in F$ nor $(X/A) \in F$. Since $S$ is an single ultrafilter, either $A \in S$ or $(X/A) \in S$. Without loss of generality, assume $A \in S$. Then, by condition (S1), for any $e \in A$ and for any $A' \in S$ such that $f(A' \cap (X/\{e\})) \leq k$, we have $A' \cap (X/\{e\}) \in S$. In particular, this means that $A \cap (X/\{e\}) \in S$ for any $e \in$ A and for any $A' \in S$ such that $f(A' \cap (X/\{e\})) \leq k$.

Now, consider the set $B = A \cup \{e\}$ for some e $\in$ A. By the submodularity of $f$, we have:

$f(A) + f(\{e\}) \geq f(A \cup \{e\}) + f(A \cap \{e\})$ .

Since $f(\{e\}) \leq k$, we have $f(A \cup \{e\}) \leq f(A) + k$. Similarly, since $f(A) \leq k$, we have $f(A \cap \{e\}) \leq k$. Therefore, $f(A \cup \{e\}) + f(A \cap \{e\}) \leq 2k + f(A)$.

Since $F$ is a single filter of order $k + 1$, it follows that $B \in F$. But then, by condition (S1) of Definition 4, we have:

$A \cap (X/\{e\}) = B \cap (X/\{e\}) \in F$

This contradicts the assumption that neither $A \in F$ nor $(X/A) \in F$, since $A \cap (X/\{e\}) \in S$ and $S \subset F$. Therefore, our assumption that such a filter F exists must be false, and we conclude that $S$ is a maximal single filter of order $k + 1$.

Therefore, we have shown that if S is an order $k + 1$ single ultrafilter, then S is a maximal single filter of order $k + 1$ under the assumption that $f(\{e\}) \leq k$ for every $e \in X$. This completes the proof.

**Lemma 16:** Let $X$ represent a finite set and $f$ denote a symmetric submodular function. Under the assumption that $f(\{e\})\leq k$ for every $e \in X$, If subset $S \subseteq 2^X$ is order $k + 1$ maximal single filter on a connectivity system $(X, f)$ , then subset $S \subseteq 2^X$ is single ultrafilter on $(X, f)$.

**Proof:** Assume that subset $S \subseteq 2^X$ is an order $k+1$ maximal single filter on a connectivity system $(X, f)$. We will prove that $S$ is a single ultrafilter on $(X, f)$.

We show that $S$ satisfies condition (S4) of Definition 5.

Let $A \subseteq X$ with $f(A) \leq k$. We need to show that either $A \in S$ or $(X/A) \in S$. Suppose that $A \notin S$. Then, by the maximality of $S$, there exists some $B \subseteq X$ such that $A \subset B$ and $B \in S$. Since $f(B) \leq k$, by condition (S2), we have $(X/B) \in S$.

We claim that $X/A = (X/B) \cap (X/\{e\})$ for some $e \in B / A$. To see this, let $e$ be any element of $B / A$. Since $A \subset B$, we have $B = A \cup (B / A)$. Thus, we have $X/B = (X/A) \cup (X/(B/A))$. Since $e \in B / A$, we have $X/B \subseteq (X/A) \cup (X/\{e\})$. On the other hand, since $A \subset B$ and $e \in B$, we have $X/\{e\} \subseteq X/(B/A)$. Therefore, we have $(X/A) \cap (X/\{e\}) = (X/B) \cap (X/\{e\})$, which implies that $X/A = (X/B) \cap (X/\{e\})$.

Now, since $B \in S$ and $f(\{e\}) \leq k$, by condition (S1) of Definition 4, we have $B \cap (X/\{e\}) \in S$. Thus, by the claim above, we have $(X/A) \in S$. Therefore, $S$ satisfies (S4).

Thus, $S$ satisfies all the conditions of Definition 5, and therefore $S$ is a single ultrafilter on $(X, f)$. This completes the proof.

An equivalent form of a given single ultrafilter $U\subseteq 2^X$ on a connectivity system is clearly manifested as a two-valued morphism. This relationship is defined through a function $m$ where $m(A) = 1$ if A is an element of $U$, and $m(A) = 0$ otherwise. Additionally, the function m holds the value $1$ for subsets where the value of the submodular function is at most $k$.

## 4. Cryptomorphism between Linear Tangle and Single Ultrafilter

In this section, we consider about cryptomorphism between linear tangle and single ultrafilter.

**Lemma 17:** Let $X$ represent a finite set and $f$ denote a symmetric submodular function. Under the assumption that $f(\{e\})\leq k$ for every $e \in X$, let $L$ be a linear tangle of order $k+1$ on $(X,f)$, and let S = {A|X/A ∈ L}. Then, S is a single ultrafilter of order k+1 on $(X,f)$.

**Proof:** Let $L$ be a linear tangle of order $k+1$, and let $S = \{A|X/A \in L\}$. We need to show that S satisfies the conditions of an order k+1 single ultrafilter on $(X,f)$.

(S1) Let $A \in S$ and $e \in X$ such that $f(\{e\}) \leq k$ and $f(A \cap (X/\{e\})) \leq k$. We need to show that $A \cap (X/\{e\}) \in S$, i.e., $X/(A \cap (X/\{e\})) \in L$. Since $A \in S$, we have $X/A \in L$. By the definition of a linear tangle, we know that there exists a partition $P$ of $X/A$ into two disjoint nonempty subsets such that for any $S,T \in P$, $X/(S \cap T) \in L$. Let $P' = \{T \cap (X/\{e\}) \mid T \in P\}$. Since $e \notin A$, we have $X/(A \cap (X/\{e\})) = X/A \cup X/(A \cap (X/\{e\}))$. Thus, $P'$ is a partition of $X/(A \cap (X/\{e\}))$ into two disjoint nonempty subsets. Moreover, for any $S,T \in P'$, we have $S \cap T = (S \cap (X/\{e\})) \cap (T \cap (X/\{e\}))$. Since $X/(S \cap (X/\{e\})) \in L$ and $X/(T \cap (X/\{e\})) \in L$, we have $X/(S \cap T) = X/((S \cap (X/\{e\})) \cup (T \cap (X/\{e\}))) = X/(S \cap (X/\{e\})) \cap X/(T \cap (X/\{e\})) \in L$. Therefore, $A \cap (X/\{e\}) \in S$.

(S2) Let $A \in S$ and $A \subset B \subseteq X$ such that $f(B) \leq k$. We need to show that $B \in S$, i.e., $X/B \in L$. Since $A \in S$, we have $X/A \in L$. Since $B$ contains $A$, we can write $B = A \cup C$ for some $C \subseteq X$. Since $f$ is submodular, we have $f(B) + f(A) \geq f(A \cup C) + f(A \cap C)$. Since $f$ is symmetric, we also have $f(B) + f(X/A) \geq f(X) + f(C)$. Adding these two inequalities, we get $f(B) + f(X) \geq f(X/A \cap C) + f(A \cup C)$. Since $A \cap C \subseteq A$ and $X/A \cap C \supseteq X/A$, we have $f(X/A \cap C) + f(A \cup C) \leq f(X) + f(A)$, which implies $f(B) + f(X) \geq f(X) + f(A)$. Thus, we have $f(B) \geq f(A)$. Since $X/A \in L$ and $f(B) \leq k$, we have $X/B = (X/A)/(B/A) \in L$. Therefore, $B \in S$.

(F3) $\emptyset$ is not in $S$ by the definition of $S$.

(S4) Let $A \subseteq X$ such that $f(A) \leq k$. We need to show that either $A \in S$ or $(X/A) \in S$. If $A \in S$, then we are done. Otherwise, we have $X/A \notin L$, which means there exists a partition $P$ of $X/A$ into two disjoint nonempty subsets such that for any $S,T \in P$, $X/(S \cap T) \notin L$. Let $P' = \{S \cup A \mid$ $S \in P\}$. Since $P$ is a partition of $X/A$, we have that $P'$ is a partition of $X$ into nonempty subsets. Moreover, for any $S,T \in P'$, we have $S \cap T = (S \cap X/A) \cap (T \cap X/A)$. Since $X/(S \cap X/A) \notin L$ and $X/(T \cap X/A) \notin L$, we have $X/(S \cap T) = X/((S \cap X/A) \cup (T \cap X/A)) = X/(S \cap X/A) \cap X/(T \cap X/A) \notin L$. Therefore, $X/A \cup A = X/(X/A) \in S$.
Thus, we have shown that $S$ satisfies all the conditions of an order $k+1$ single ultrafilter on $(X,f)$. This completes the proof.

**Lemma 18:** Let $X$ represent a finite set and $f$ denote a symmetric submodular function. Under the assumption that $f(\{e\})\leq k$ for every $e \in X$, let $S$ be a single ultrafilter of order $k+1$, and let $L = \{A|\ (X/A) \in S\}$. Then, $L$ is a linear tangle of order $k+1$.
**Proof:** Let $S$ be an order $k + 1$ single ultrafilter of $X$ and let $L = \{A \mid (X/A) \in S\}$. We need to show that $L$ is a linear tangle of order $k + 1$.
First, we show that $\varnothing$ belongs to $L$. Since $S$ is an ultrafilter, $\varnothing$ is not in $S$, which implies that $X/\varnothing = X$ is in $S$. Thus, $\varnothing$ belongs to $L$.
Next, we show that for any $k$-efficient subset $A$ of $X$, either $A$ or $X/A$ belongs to $L$. Suppose $A$ is a $k$-efficient subset of $X$. If $A$ is not in $L$, then $X/A$ is in $S$. By Definition 5, we have that for any $e \in X$, if $f(\{e\}) \leq k$ and $f(A \cap (X/\{e\})) \leq k$, then $A \cap (X/\{e\})$ belongs to $S$. Since $A$ is $k$-efficient, $f(A) \leq k$. Therefore, for any $e \in X, f(\{e\}) \leq k$ and $f(A) \leq k,$ which implies that $f(A \cap (X/\{e\})) \leq k$. Thus, by Definition 5, $A \cap (X/\{e\})$ belongs to $S$ for any $e \in A$. Since $A$ is a $k$-efficient subset of $X$, $f(X/A) = f(A) \leq k$. Hence, by Definition 5, $X/A$ belongs to $S$. Therefore, either $A$ or $X/A$ belongs to $L$.
Lastly, we show that for any $A, B \in L$ and $e \in X$ such that $f(\{e\}) \leq k$, the set $A \cup B \cup \{e\}$ is not equal to $X$. Let $A$, B be two elements of $L$, and let $e$ be an element of $X$ such that $f(\{e\}) \leq k$. Suppose $A \cup B \cup \{e\} = X$. Then, we have that $(X/(A \cup B)) \subseteq (X/\{e\})$. Since $A, B$ belong to $L$, we have that $(X/A), (X/B)$ belong to $S$. Therefore, we have that $(X/(A \cup B)) = (X/A) \cap (X/B)$ belongs to $S$ by Definition 5. Hence, $(X/\{e\})$ contains an element that belongs to $S$, which contradicts the fact that $f(\{e\}) \leq k$ and $S$ is an order $k + 1$ single ultrafilter of $X$. Therefore, $A \cup B \cup \{e\}$ is not equal to $X$.
Thus, we have shown that $L$ satisfies all the axioms of a linear tangle of order $k + 1$. Therefore, $L$ is a linear tangle of order $k + 1$.This completes the proof.

From above lemmas, we obtain the following:
**Theorem 19:** Let $X$ represent a finite set and $f$ denote a symmetric submodular function. Under the assumption that $f(\{e\})\leq k$ for every $e \in X$, the family $L$ being a $k+1$ order single ultrafilter and $L$ being a $k+1$ order linear tangle are equivalent necessary and sufficient conditions.

Linear tangle are known as obstructions to linear decomposition [1]. Therefore, Theorem 6 shows that single ultrafilter are obstructions to linear decomposition under the assumption that $f(\{e\})\leq k$ for every $e \in X$.

## 5. Matroid, Antimatroid, Greedoid on a connectivity system *(X,f)*

Matroids are a discrete representation of linear algebra concepts such as linearly independent and linearly dependent sets. They are fundamental and significant concepts closely related to computer science. As a result, there has been a wealth of literature published on the topic [42-47]. The concepts of matroid theory are well-suited for greedy algorithms, making them a useful tool in various algorithmic applications. This underscores their importance and relevance. Additionally, related concepts such as Antimatroids and Greedoids are also well-known and actively researched.

In literature[52], we present a definition of matroids on connectivity systems $(X, f)$ that mimics the concept of matroids on Boolean algebra. By imposing restrictions on submodular functions for each axiom, we ensure a coherent and natural definition.

**Definition 21:** In a connectivity system $(X,f)$, the set family $M \subseteq 2^X$ is called a matroid of order $k+1$ *on* $(X,f)$ if the following axioms hold true:
(M0) For every $A \in M, f(A) \leq k,$

(M1) $\varnothing \in M$,
(M2) if $A \in M, B \subseteq A$, and $f(B) \leq k$ then $B \in M$,.
(M3) if $A, B \in M, |A| < |B|, e \in X, f(\{e\}) \leq k$, and f(A $\cup$ {e}) $\leq k$, , then $e \in B / A$ such that $A \cup \{e\} \in M$.

Furthermore, we introduce new concepts beyond the ones mentioned above. In this paper, we extend the notion of Greedoid from matroid theory to the connectivity system *(X, f)* by defining that a subset $M \subseteq 2^X$ is a Greedoid on *(X, f)* if it satisfies properties (M0), (M1), and (M3).

Additionally, an Antimatroid on the connectivity system *(X, f)* is defined as a subset $M \subseteq 2^X$ that satisfies properties (M0), (M1), (M3), and (AM2). By naturally extending the concept of Antimatroid from matroid theory, we incorporate it into the connectivity system *(X, f)*.
(AM2) If $A \in M$, $A \neq \emptyset$, $e \in X$, $f(\{e\}) \leq k$, $f(A \setminus \{e\}) \leq k$ then $e \in A$ such that $A \setminus \{e\} \in M$.

Let us define an order *k+1* matroid (antimatroid, greedoid respectively) $M \subseteq 2^X$ on connectivity system *(X,f)* as an Ultra Matroid (Ultra antimatroid, Ultra greedoid respectively) if it satisfies the following Axiom (M4):
(M4) For any subset $A \subseteq X$, if $f(A) \leq k$, then either $A \in M$ or $X/A \in M$.

Axiom (M4) imitates the idea of Ultrafilters, which states that "either a set A or its complement belongs to the (Ultra)filter." Therefore, ultra matroids (Ultra antimatroid, Ultra greedoid respectively) can be interpreted as matroids (antimatroids, greedoids respectively) that include only one side.

Additionally, we introduce a new concept. A matroid (Ultra Matroid) on a connectivity system *(X, f)* is considered prime if it satisfies the following axiom:

$$\text{(M5) } A \in M \text{ for all } A \subseteq X \text{ with } |A| = 1, f(A) \leqq k.$$

**Theorem 22:** Let *X* be a finite set and *f* be a symmetric submodular function. Linear-width of the connectivity system *(X, f)* is at most *k* if and only if no prime Ultra matroid of order *k+1* exists.

**Proof:** This proof is based on techniques from the paper [12]. The context begins with a connectivity system *(X, f)* and a prime Ultra matroid *M* of order *k+1*. A subset A of *X* is said to be "*k-l*inear-branched" if the connectivity system resulting from collapsing *X\A* has a linear width of at most *k*. (If *X* itself is *k*-linear-branched, it is equivalent to the linear width of *f* being at most *k*.) *B* represents the set of all *k*-branched subsets of *X. B'* is then defined as the set of all *A* such that *A* is a subset of some set *C* in *B,* with the constraint that *f(A)* is less than or equal to *k*.

Now let's consider the two directions needed to prove this theorem:

**Direction 1:** If no prime Ultra matroid of order *k+1* exists, the Linear-width of the connectivity system *(X, f)* is at most *k*. Here, we're looking at the condition where there is no prime Ultra matroid of order *k+1*. If this condition holds true, we can conclude that the linear-width of the connectivity system *(X, f)* is at most *k*. This is the "if" direction of the proof.
**Direction 2:** Conversely, if the linear width of the connectivity system *(X, f)* is at most *k,* then *X* is *k*-linear-branched, and therefore, *X* belongs to *B'*.

In this case, we're starting with the assumption that the linear width of the connectivity system *(X, f)* is at most *k*. If this is the case, then we know that *X* is *k*-linear-branched. This in turn means that *X* belongs to the set *B'*. This is the "only if" direction of the proof.

Moreover, it is clear that every *k*-linear-branched set with at least two elements can be expressed as the union of two proper *k*-linear-branched subsets. This means that no prime Ultra matroid of order *k+1* exists.
By proving these two directions, we have shown that the linear width of the connectivity system *(X, f)* being at most *k* is equivalent to there being no prime Ultra matroid of order *k+1*. This completes the proof of Theorem 22.

**Theorem 23:** Let *X* be a finite set and *f* be a symmetric submodular function. Linear-width of the connectivity system *(X, f)* is at most *k* if and only if no prime Ultra antimatroid of order *k+1* exists.
**Proof:** The proof follows a similar approach as in Theorem 22.

**Theorem 24:** Let *X* be a finite set and *f* be a symmetric submodular function. Linear-width of the connectivity system *(X, f)* is at most $k$ if and only if no prime Ultra greedoid of order $k+1$ exists.

**Proof:** The proof follows a similar approach as in Theorem 22.
We consider Quasi-Matroids on a connectivity system (X,f) [53]. Quasi-Matroids, as defined in reference [54] (cf. [55]), represent a conceptual extension of Matroids, known for their generalized nature. When extended to the realm of a connectivity system, the following holds true:

**Theorem 25:** Let $X$ be a finite set and f a symmetric submodular function. The linear-width of the connectivity system $(X, f)$ is at most $k$ if and only if no prime Ultra Quasi-Matroid of order $k+1$ exists.
**Proof:** The proof follows a similar approach as in Theorem 22.

## 6. Conclusion: Duality theorem

From the previous results in this paper, the following duality theorem can be obtained.

**Theorem 26:** Let $X$ be a finite set and $f$ be a symmetric submodular function. Linear-width of the connectivity system $(X, f)$ is at most $k$ if and only if following conditions:
1. No prime Ultra antimatroid of order $k+1$ on connectivity system $(X, f)$ exists.
2. No prime Ultra matroid of order $k+1$ on connectivity system $(X, f)$ exists.
3. No prime Ultra greedoid of order $k+1$ on connectivity system $(X, f)$ exists.
4. No Linear tangle of order $k+1$ on connectivity system $(X, f)$ exists.
5. No Linear obstacle of order $k+1$ on connectivity system $(X, f)$ exists.
6. No maximal single ideal of order $k+1$ on connectivity system $(X, f)$ exists.
7. No Maximal Linear loose tangle of order $k+1$ on connectivity system $(X, f)$ exists.
8. No (non-principal) single ultrafilter of order $k+1$ on connectivity system $(X, f)$ exists.
9. No Ultra Quasi-matroids of order $k+1$ on connectivity system $(X, f)$ exists.

## 7. Future Tasks: Superfilter

In the future, we plan to investigate the relationship between the following definition, which adds the condition of symmetric submodular function to the definition on Boolean algebra, and the graph width parameter to further deepen our understanding of the graph width parameter. This definition is well-known as a concept closely related to filters in Boolean algebra.
The definitions provided above serve as mere examples of concepts, and there exist numerous other concepts that can be extended and examined within the framework of a connectivity system.

**Definition 27:** Let $X$ be a finite set and $f$ be a symmetric submodular function. A set family $F \subseteq 2^X$ is called an order $k + 1$ **closed set system** on a connectivity system $(X, f)$ if it satisfies the following conditions:
(C1) $\emptyset$ is not belong to $F$,
(C2) $A, B \in F, f(A \cap B) \leq k \Rightarrow A \cap B \in F$,
(C3) $A \subseteq B \subseteq X, f(A) \leq k, A \cup B \in F, f(B) \leq k \Rightarrow B \in F$.

**Definition 28(cf.[60]):** Let $X$ be a finite set and $f$ be a symmetric submodular function. A set family $F \subseteq 2^X$ is called an order $k + 1$ **superfilter** on a connectivity system $(X, f)$ if it satisfies the following conditions:
(S1) If $A \in S, B \supseteq A$, and $f(B) \leq k$ then $B \in S$.
(S2) If $A \cup B \in S, f(A) \leq k$ *and* $f(B) \leq k$ then $A \in S$ or $B \in S$.

Also, we aim to further advance the study of graph parameters by adding the submodularity condition to various concepts on algebra, such as Frechet filter, Complemented filter, Upper set, lower set, Scott set, filter base, monotone class, subalgebra, prime filter, complemented ideal, matroid (cf.[58]), and others.

## Acknowledgments

I humbly express my sincere gratitude to all those who have extended their invaluable support, enabling me to successfully accomplish this paper.